\documentclass[a4paper,12pt]{article}

\usepackage{verbatim}
\usepackage{amssymb}
\usepackage{amsbsy}
\usepackage{amscd}
\usepackage{amsmath}
\usepackage{amsthm}
\usepackage[mathscr]{eucal}
\usepackage{cite}

\usepackage{color}
\usepackage{latexsym}
\usepackage{url}
\usepackage{multirow}
\usepackage{booktabs}
\usepackage{enumitem}

\def\d{\displaystyle}
\def\p{\partial}
\def\s{\sigma}
\def\e{\varepsilon}

\def\R{{\bf R}}
\def\N{{\bf N}}

\newcommand{\norm}[1]{\left\lVert#1\right\rVert}
\DeclareMathOperator*{\supp}{supp}

\newtheorem{theorem}{Theorem}
\newtheorem{definition}{Definition}\numberwithin{definition}{section}
\newtheorem{lemma}{Lemma}\numberwithin{lemma}{section}
\newtheorem{remark}{Remark}\numberwithin{remark}{section}

\title{Short time blow-up by negative mass term for semilinear wave equations with small data and scattering damping}

\author{Ning-An Lai
\footnote{Institute of Nonlinear Analysis and Department of Mathematics, Lishui University, Lishui City
	323000, China. e-mail : ninganlai@lsu.edu.cn.},
Nico Michele Schiavone
\footnote{Department of Mathematics, University of Rome ``La Sapienza'', Piazzale Aldo Moro 5, 00185 Roma, Italy. e-mail : schiavone@mat.uniroma1.it.},
Hiroyuki Takamura
\footnote{Mathematical Institute, Tohoku University, Aoba, Sendai 980-8578, Japan.
e-mail : hiroyuki.takamura.a1@tohoku.ac.jp.}}

\date{
\footnotesize
Keywords : wave equation, semilinear, damping, mass, blow-up, lifespan\\
MSC2010 : primary 35L71, secondary 35B44}

\begin{document}

\maketitle

\begin{abstract}
In this paper we study blow-up and lifespan estimate for solutions to the Cauchy problem with small data for semilinear wave equations with scattering damping and negative mass term.
We show that the negative mass term will play a dominant role when the decay of its coefficients is not so fast, thus the solutions will blow up in a finite time. What is more, we establish a lifespan estimate from above which is much shorter than the usual one.
\end{abstract}

\section{Introduction}

We consider the Cauchy problem with small data for the semilinear wave equations with scattering damping and negative mass term
\begin{equation}
\label{eq:main_problem}
\left\{
\begin{aligned}
& u_{tt} - \Delta u + \frac{\mu_1}{(1+t)^\beta} u_t - \frac{\mu_2}{(1+t)^{\alpha+1}} u = |u|^p, \quad\text{in $\R^n\times(0,T)$}, \\
& u(x,0)=\e f(x), \quad u_t(x,0)=\e g(x), \quad x\in\R^n,
\end{aligned}
\right.
\end{equation}
where  $\mu_1\ge0$, $\mu_2 > 0$, $\alpha<1$, $\beta >1$, $p>1$, $n\in\N$, $T>0$ and $\e>0$ is a ``small'' parameter. This problem comes from the recent interest in the \lq\lq wave-like" or \lq\lq heat-like" behaviour of semilinear wave equations with variable coefficients
damping. The Cauchy problems with small data for
\begin{equation*}
	u_{tt}-\Delta u=|u|^p
	\quad\text{and}\quad
	u_t-\Delta u=|u|^p
\end{equation*}
admit critical powers, respectively, the so-called Strauss exponent $p_S(n)$ and the Fujita exponent $p_F(n)$ (see \cite{Strauss} and \cite{Fujita}), where for \lq\lq critical power" of a problem we mean the exponent $p_{c}$ such that its small data solutions blow up for $1<p\le p_{c}$ and exist globally in time for $p>p_{c}$.
% means that this exponent divides the region $1<p<\infty$ into two parts, one (smaller than it) corresponds to the blow-up of solutions and the other one corresponds to global existence.
It is of recent interest the Cauchy problem with small data for
\begin{equation}
\label{dampedeq}
\begin{split}
u_{tt}-\Delta u+\frac{\mu}{(1+t)^\beta}u_t=|u|^p.
\end{split}
\end{equation}
If the Cauchy problem \eqref{dampedeq} admits a critical power related to $p_S(n)$, then we say it has a \lq\lq wave-like'' behaviour, while if it is related to $p_F(n)$, then we say it admits a \lq\lq heat-like'' behaviour.
Generally speaking, if the decay rate $\beta$ of the damping coefficients is large enough, then the damping term seems to have no influence and then we get a \lq\lq wave-like'' behaviour; otherwise, we get a \lq\lq heat-like'' behaviour.
According to the different value of $\beta$, we recover four cases  (overdamping, effective, scaling invariant, scattering), based on the works by Wirth \cite{Wir1, Wir2, Wir3} (see also \cite{DLR15, FIW, II, IWnew, LTW, LT, TL1, WY, WY14_scale, Wang} and references therein).

\begin{comment}
In this sense, there exists also a critical decay rate of the variable coefficients to distinguish between wave and heat. The known results (see \cite{DLR15}, \cite{FIW}, \cite{II}, \cite{IWnew}, \cite{LTW}, \cite{LT}, \cite{TL1}, \cite{WY}, \cite{WY14_scale}, \cite{Wang} and references therein) indicate that the threshold is $\beta=1$. Actually, we have four cases according to the different value of $\beta$ (overdamping, effective, scaling invariant, scattering), based on the works by Wirth \cite{Wir1, Wir2, Wir3} on the corresponding linear problem for
	\begin{equation*}
	u^0_{tt}-\Delta u^0+\frac{\mu}{(1+t)^\beta}u^0_t=0.
	\end{equation*}
\end{comment}

\begin{comment}
\begin{equation*}
\label{linearproblem}
\left\{
\begin{array}{l}
\d u^0_{tt}-\Delta u^0+\frac{\mu}{(1+t)^\beta}u^0_t=0, \quad \mbox{in}\ \R^n \times[0,\infty), \\
u^0(x,0)=u_1(x),\ u^0_t(x,0)=u_2(x), \quad x\in\R^n.
\end{array}
\right.
\end{equation*}
\end{comment}
\begin{comment}
We list the four cases in the following table.
\begin{center}
	\begin{tabular}{|c|c|}
		\hline
		$\beta\in(-\infty,-1)$ & overdamping\\
		\hline
		$\beta\in[-1,1)$ & effective\\
		\hline
		$\beta=1$ &
		\begin{tabular}{c}
			scaling invariant\\
			$\mu\in(0,1)\Rightarrow$ non-effective
		\end{tabular}
		\\
		\hline
		$\beta\in(1,\infty)$ & scattering\\
		\hline
	\end{tabular}
\end{center}
\end{comment}

\par
On the other hand, people are paying more attention to the Cauchy problem for
\begin{equation*}
	u_{tt} - \Delta u + \frac{\mu_1}{1+t} u_t + \frac{\mu_2^2}{(1+t)^{2}} u = |u|^p,
\end{equation*}
\begin{comment}
\begin{equation*}
\label{nonlinearmass}
\left\{
\begin{aligned}
& u_{tt} - \Delta u + \frac{\mu_1}{1+t} u_t + \frac{\mu_2^2}{(1+t)^{2}} u = |u|^p,
\quad \text{in $\R^n \times[0,\infty)$},\\
& u(x,0)=u_1(x), \quad u_t(x,0)=u_2(x), \quad x\in\R^n,
\end{aligned}
\right.
\end{equation*}
\end{comment}
which includes scale-invariant damping and mass in the mean time. In some sense, this model describes the interplay between the damping and mass. For this problem, the quantity
\begin{equation*}
\label{eq:def_delta}
\delta:=(\mu_1-1)^2-4\mu_2^2
\end{equation*}
plays an important role to the behaviour of the solutions. We refer the reader to \cite{NPR, Pal, P-R} and references therein.

Naturally, we want to consider the corresponding problem with scattering damping and mass term. Very recently, the authors \cite{LST}
studied the Cauchy problem \eqref{eq:main_problem} with fast decay rate in the coefficients of the mass term, thus, $\alpha>1$, in which we proved
\begin{comment}
blow-up results for
\begin{equation*}
	p>1,\,n=1 \quad\text{and}\quad 1<p<p_S(n),\, n\ge 2.
\end{equation*}
\end{comment}
that the blow-up results and the upper bound lifespan estimates are the same as that of the semilinear wave equations with scattering damping but without mass term, see \cite{LT}. This implies that the negative mass term seems to have no influence on the behaviour of the solutions.
In this work, we are devoted to studying the case $\alpha<1$.
Our motivation to study a negative mass term is simply a mathematical interest by \cite{LST},
but one may refer to the introduction of Yagdjian and Galstian \cite{YG} which mentions its physical background.
From our main result listed in Theorem \ref{thm:alpha<1} below, it seems that the negative mass term will affect the qualitative properties of the small data solutions of the Cauchy problem \eqref{eq:main_problem}, due to two reasons:
firstly, the non-existence of global energy solutions can be established for all $p>1$ and $n\ge 1$; moreover, the upper bound of the lifespan is smaller than the usual one and it looks like a log-type with respect to $\e$.
%\textcolor{red}{We can also observe that the scattering damping plays no role (we will absorb it by a multiplier, first introduced in \cite{LT}).}
%\textcolor{red}{Observe that if we set in particular $\mu_1=0$, $\alpha=-1$ and $m=\sqrt{\mu_2}$, our equation reduces to $u_{tt}-\Delta u + (im)^2 u=|u|^p$, that is the Klein-Gordon equation with power nonlinearity for imaginary mass particles (tachyons).}

\begin{comment}
	In order to prove our main result, we use a multiplier, which was first introduced in \cite{LT}, to absorb the damping term. Then we establish the lifespan from above by improved Kato's type lemma, following the idea in \cite{T}.
\end{comment}

%%%%%%%%%%%%%%%%%%%
%%%%%% section 2 %%%%%%%
%%%%%%%%%%%%%%%%%%%

\section{Definitions and Main Result}
	First of all, let us introduce energy
	%and weak
	solutions of Cauchy problem \eqref{eq:main_problem}.

\begin{definition}
	We say that $u$ is an energy solution for problem \eqref{eq:main_problem} over $\R^n\times[0,T)$ if
	\begin{equation*}
	u \in C([0,T), H^1(\R^n)) \cap C^1([0,T),L^2(\R^n)) \cap L^p_{loc}(\R^n \times (0,T))
	\end{equation*}
	satisfies $u(x,0)=\e f(x)$ in $H^1(\R^n)$, $u_t(x,0) = \e g(x)$ in $L^2(\R^n)$ and
	\begin{equation}
	\label{eq:energy_solution}
	\begin{split}
	&\int_{\R^n}u_t(x,t)\phi(x,t)dx-\int_{\R^n}\e g(x)\phi(x,0)dx\\
	&+\int_0^tds\int_{\R^n}\left\{-u_t(x,s)\phi_t(x,s)+\nabla u(x,s)\cdot\nabla\phi(x,s)\right\}dx \\
	&+\int_0^t ds \int_{\R^n}
	\left\{ \frac{\mu_1 }{(1+s)^{\beta}}u_t(x,s)
	- \frac{\mu_2}{(1+s)^{\alpha+1}}u(x,s)
	\right\}\phi(x,s)dx
	\\
	= & \int_0^tds\int_{\R^n}|u(x,s)|^p\phi(x,s)dx
	\end{split}
	\end{equation}
	for all test functions $\phi \in C_0^\infty(\R^n \times [0,T))$ and for all $t\in[0,T)$.
\end{definition}

\begin{comment}
Employing the integration by parts in the above equality and letting $t\to T$, we got the definition of the weak solution of \eqref{eq:main_problem}, that is
\begin{equation*}
\begin{split}
&\int_{\R^n\times[0,T)}
u(x,s) \Bigg\{\phi_{tt}(x,s)-\Delta \phi(x,s) \\
&-\frac{\p}{\p s} \left(\frac{\mu_1}{(1+s)^{\beta}}\phi(x,s)\Bigg)
- \frac{\mu_2}{(1+s)^{\alpha+1}} \phi(x,s) \right\}dxds\\
=&\ \int_{\R^n}\mu_1 \e f(x)\phi(x,0)dx -\int_{\R^n}\e f(x)\phi_t(x,0)dx +\int_{\R^n}\e g(x)\phi(x,0)dx\\ &+\int_{\R^n\times[0,T)}|u(x,s)|^p\phi(x,s)dxds.
\end{split}
\end{equation*}
\end{comment}

%%%%%%%%%%%%

\begin{theorem}\label{thm:alpha<1}
	Let $\alpha<1$, $n\ge1$ and $p>1$. Assume that both $f \in H^1(\R^n)$ and $g\in L^2(\R^n)$ are non-negative, at least one of them does not vanish identically.
		 Suppose that $u$ is an energy solution of \eqref{eq:main_problem} on $\R^n \times [0,T)$ that satisfies, for some $R\ge1$,
	\begin{equation}
	\label{support}
	\supp u \subset\{(x,t)\in\R^n\times[0,T) \colon |x|\le t+R\}.
	\end{equation}

	Then, there exists a constant $\e_0=\e_0(f,g,R,n,p,\mu_1,\mu_2,\alpha,\beta)>0$
	such that $T$ has to satisfy
	\begin{equation*}
	%\label{Testweak}
	T\leq 3 \zeta(C\e)
	\end{equation*}
	for $0<\e\le\e_0$, where $\zeta=\zeta(\overline{\e})$ is the larger solution to the equation
	\begin{equation}
	\label{eq:defb}
	\overline{\e} \zeta^{\frac{2}{p-1}-n+\frac{1+\alpha}{4}}
	\exp\left(\frac{2}{1-\alpha}\sqrt{ \mu_2 \exp\left(\frac{\mu_1}{1-\beta}\right) } \zeta^{\frac{1-\alpha}{2}}\right) = 1
	\end{equation}
	and $C$ is a positive constant independent of $\e$.
\end{theorem}

\begin{remark}\label{rem1}
Let us make some observations:
\begin{itemize}
\item The assumption (\ref{support}) can be replaced by
$\supp\{f,g\} \subset \{x\in\R^n \colon |x|\le R \}$ when $n=1,2$, or $p\le n/(n-2)$ for $n\ge3$.
This fact is established by local existence of such an energy solution.
See Appendix in the last section.
\item
Since letting $\e\to0$ we have $\zeta\to+\infty$ in (\ref{eq:defb}),
 it is not difficult to see that $T \le c [\log(1/\e)]^{2/(1-\alpha)}$ for some constant $c>0$
follows from (\ref{eq:defb}).
In fact, this inequality is trivial when the exponent $\delta:=2/(p-1)-n+(1+\alpha)/4$ of the first $\zeta$
is non-negative,
while $\zeta^\delta$ can be absorbed by square root of the exponential term when $\delta<0$.
\item
It is an open question the optimality of the upper bound of the lifespan in Theorem~\ref{thm:alpha<1}.
\end{itemize}
\end{remark}

%removed from the previous remark
\begin{comment}
To this end, one may rewrite the equation as
\[
u_{tt} - \Delta u + \frac{\mu_2}{(1+t)^{\alpha+1}} u = |u|^p-\frac{\mu_1}{(1+t)^\beta} u_t+\frac{2\mu_2}{(1+t)^{\alpha+1}} u
\]
\end{comment}

%%%%%%%%%%%%%%%%%%%
%%%%%% section 3 %%%%%%%
%%%%%%%%%%%%%%%%%%%

\section{Kato's type lemma}
\label{sec:kato}
In order to prove our theorem, we need a slightly different version of the improved Kato's lemma introduced in \cite{T}.

\begin{lemma}
	\label{lem:kato}
	Let $p>1$ and $0\le \widetilde{T}_0<T$ be positive constants. Suppose that $A \in C^1([\widetilde{T}_0,T))$, $B \in C^1([0,T))$, $m \in C^1([0,T))$ are strictly positive functions, $B=B(t)$ is decreasing and that $m(t)$ is bounded by two constants $\overline{m}, \underline{m}>0$($\underline{m} \le m(t) \le \overline{m}$) for $t\ge0$.
	Define the function
	\begin{equation}
	\label{def:h}
	h(t) := B(t)^{1/2} A(t)^{(p-1)/2-\delta}, %\quad
	%h_0(t) := B(t)^{1/2} A(t)^{(p-1)/2},
	\end{equation}
	where $\delta$ is a constant such that $0<\delta<(p-1)/2$ and $h'(t) \ge 0$ for $t \ge \widetilde{T}_0$.
	
	Assume that $F \in C^2([0,T))$ satisfies
	\begin{gather}
	\label{hp1}
	F(t) \ge A(t) \quad\text{for $t\ge \widetilde{T}_0$}, \\
	\label{hp2}
	\{m(t)F'(t)\}' \ge B(t) |F(t)|^p \quad\text{for $t\ge 0$}, \\
	\label{hp3}
	F(0), F'(0) \ge 0, \quad  F(0) + F'(0)>0.
	\end{gather}
	If $F'(0)=0$, suppose that there exists a time $\widetilde{t} > 0$ such that
	\begin{equation}
	\label{F(t)>2F0}
	F(\widetilde{t}) \ge 2 F(0) .
	\end{equation}
	Define the time %$\widetilde{T}_1:=\overline{m}\underline{m}^{-1} {F(0)}/{F'(0)}$ if $F'(0) \neq 0$, and $\widetilde{T}_1:=\widetilde{t}$ if $F'(0)=0$.
	
		\begin{equation*}
		\label{def:wT1}
		\widetilde{T}_1:=
		\begin{cases}
		\overline{m}\underline{m}^{-1} {F(0)}/{F'(0)} & \text{if $F'(0) \neq 0$,} \\
		\widetilde{t} & \text{if $F'(0)=0$.}
		\end{cases}
		\end{equation*}

	Then,
	for $\widetilde{T} \ge \max\{\widetilde{T}_0,\widetilde{T}_1\}$ we have
	%\begin{equation*}
	%\label{eq:T}
	$T \le 3\widetilde{T}$
	%\end{equation*}
	assuming that
	\begin{equation}
	\label{eq:T1estimate}
	\widetilde{T} \, h(\widetilde{T}) A(\widetilde{T})^\delta
	\ge %\frac{\overline{m}}{\delta}\sqrt{\frac{p+1}{\underline{m}}}.
	{\delta}^{-1}{\overline{m}}\sqrt{(p+1)/{\underline{m}}}.
	%\widetilde{T} \, h_0(\widetilde{T}) \ge \frac{\overline{m}}{\delta}\sqrt{\frac{p+1}{\underline{m}}}.
	\end{equation}
\end{lemma}

%\begin{remark}
%	Comparing this lemma with Lemma~\ref{lemmaT}, we observe that here we generalize the functions estimating $F$ and $F''$, but we require stronger condition on them. In fact, choosing $A(t)=At^a$, $B(t)=B(t+R)^{-b}$, the condition $h'(t)\ge0$ corresponds to $M \ge 1$.
%\end{remark}

\begin{comment}
\begin{remark}
	%\label{rmk:F>2F0}
	Before starting with the proof, observe that if $G \in C^1([t_0,t_1))$ is a positive convex function such that $G'(t_0)>0$, then
	\begin{equation*}
	t \ge t_0 + \frac{G(t_0)}{G'(t_0)} \Longrightarrow G(t) \ge 2G(t_0).
	\end{equation*}
	In fact, by the convexity, the function $G$ is above its tangent in the point $t=t_0$, i.e. $G(t) \ge G'(t_0)\{t-t_0\}+G(t_0)$, from which the conclusion.
\end{remark}
\end{comment}

\par
\begin{proof}\renewcommand{\qedsymbol}{}
	First of all, let us prove that $F(t), F'(t)>0$ for $t>0$. We need to consider two cases according to the initial conditions \eqref{hp3} on $F$.
\vskip5pt
\par\noindent
\textit{Case 1: $F'(0)>0$.} From \eqref{hp2} it follows $F'(t)\ge m(0)F'(0) m(t)^{-1} >0$, and then $F(t)\ge F(0) + m(0)F'(0) \int_{0}^{t}m(s)^{-1}ds >0$ for $t>0$.
\vskip5pt
\par\noindent
\textit{Case 2: $F'(0)=0$.} Then $F(0)>0$. It follows from \eqref{hp2} evaluated in $t=0$ that $\{mF'\}'(0)>0$, which implies $m(t)F'(t)>m(0)F'(0)=0$ for small $t>0$. Hence, the fact that $\{mF'\}'(t)\ge0$ for $t\ge0$ from \eqref{hp2} yields that $m(t)F'(t)>0$, that is $F'(t)>0$, and so $F(t)>F(0)>0$ for $t>0$.
\vskip5pt
	Moreover, observe that
	\begin{equation}
	\label{F>2F0}
	F(t) \ge 2 F(0) \quad\text{for $t\ge \widetilde{T}_1$}.
	\end{equation}
	Indeed, if $F'(0)=0$, it follows by the hypothesis \eqref{F(t)>2F0} and by the fact that $F$ is increasing. If $F'(0)>0$, by \eqref{hp2} and because $m$ is bounded, we have $F(t) \ge F(0) + \underline{m}\overline{m}^{-1} F'(0) t$, from which \eqref{F>2F0} follows.
	
	Multiplying \eqref{hp2} by $m(t)F'(t)>0$, we get
	\begin{equation*}
	\left(\frac{1}{2} \{m(t)F'(t)\}^2\right)' \ge m(t)B(t)|F(t)|^pF'(t) \quad\text{for $t>0$}.
	\end{equation*}
	From this inequality, the positivity of $F$ and the facts that $B$ is decreasing and $m$ is bounded, it follows that
	\begin{equation*}
	\begin{split}
	\frac{1}{2} F'(t)^2 &\ge
	\frac{1}{2} \overline{m}^{-2} \underline{m}^2 F'(0)^2 + \overline{m}^{-2}\underline{m} B(t) \int_{0}^{t} F(s)^p F'(s)ds \\
	%& \ge \frac{\overline{m}^{-2}\underline{m}}{p+1} B(t) \{F(t)^{p+1}-F(0)^{p+1}\} \\
	& \ge \frac{\overline{m}^{-2}\underline{m}}{p+1} B(t) F(t)^p \{F(t)-F(0)\} \quad\text{for $t\ge 0$},
	\end{split}
	\end{equation*}
	and so, using equation \eqref{F>2F0}, we get
	\begin{equation}
	\label{5}
	F'(t) \ge \overline{m}^{-1} \sqrt{{\underline{m}}/{(p+1)}} B(t)^{1/2} F(t)^{(p+1)/2} \quad \text{for $t\ge \widetilde{T}_1$}.
	\end{equation}
	
	Now, fix $\widetilde{T} := \max \{\widetilde{T}_0,\widetilde{T}_1\}$ and define the function
	\begin{equation*}
	H(t) := \int_{\widetilde{T}}^{t} h(s)ds
	= \int_{\widetilde{T}}^{t} B(s)^{1/2} A(s)^{(p-1)/2-\delta} ds
	\quad\text{for $t \ge \widetilde{T}$.}
	\end{equation*}
	%Note that from the hypothesis $h' \ge0$ we obtain that $h$ is increasing and $H$ is convex.
	Because $0<\delta<(p-1)/2$, from inequality \eqref{5} and from \eqref{hp1} we obtain
	\begin{equation*}
	\begin{split}
	{F'(t)}/{F(t)^{1+\delta}}
	%& \ge \frac{1}{\overline{m}}\sqrt{\frac{\underline{m}}{p+1}}
	%B(t)^{1/2} F(t)^{(p-1)/2-\delta} \\
	& \geq \overline{m}^{-1} \sqrt{{\underline{m}}/{(p+1)}}
	B(t)^{1/2} A(t)^{(p-1)/2-\delta} \quad \text{for $t\ge \widetilde{T}$.}
	\end{split}
	\end{equation*}
	
	Integrating this inequality on $[2\widetilde{T},t]$ we get
	\begin{equation*}
	\frac{1}{\delta} \left(\frac{1}{F(2\widetilde{T})^\delta}-\frac{1}{F(t)^\delta}\right)
	\ge \frac{1}{\overline{m}}\sqrt{\frac{\underline{m}}{p+1}}
	\int_{2\widetilde{T}}^t B(s)^{1/2} A(s)^{(p-1)/2-\delta} ds.
	\end{equation*}
	Neglecting the second term on the left-hand side, from \eqref{hp1} evaluated in $t=2\widetilde{T}$ and recalling the definition of $H$, we have
	\begin{equation}
	\label{4}
	A(2\widetilde{T})^{-\delta} \ge F(2\widetilde{T})^{-\delta}
	\ge
	\delta \overline{m}^{-1}
	\sqrt{{\underline{m}}/{(p+1)}}
	[H(t)-H(2\widetilde{T})]
	%\quad\text{for $t\ge 2\widetilde{T}$.}
	\end{equation}
	for $t\ge 2\widetilde{T}$. Since $h$ is increasing, we get
	\begin{equation}
	\label{6}
	H(2\widetilde{T})
	\ge
	h(\widetilde{T}) \int_{\widetilde{T}}^{2\widetilde{T}}ds=
	\widetilde{T}h(\widetilde{T}).
	\end{equation}
\par
	Observe moreover that $A$ is increasing, in fact
	\begin{equation*}
	h'(t) = h(t) \left\{ {B'(t)}/{(2B(t))} + \left[{(p-1)}/{2}-\delta\right] {A'(t)}/{A(t)} \right\},
	\end{equation*}
	and so, because $h,A,B >0$, $h' \ge 0$ and $B' \le 0$, we get $A' \ge 0$.
	Then, by equation \eqref{6}, hypothesis \eqref{eq:T1estimate} and the monotonicity of $A$, we have
	\begin{equation*}
	A(2\widetilde{T})^{\delta} H(2\widetilde{T}) \ge
	A(\widetilde{T})^{\delta} \, \widetilde{T} \, h(\widetilde{T})
	%= \widetilde{T} \, h_0(\widetilde{T})
	\ge %\frac{\overline{m}}{\delta}\sqrt{\frac{p+1}{\underline{m}}}
	{\delta}^{-1}{\overline{m}}\sqrt{(p+1)/{\underline{m}}}
	.
	\end{equation*}
	Inserting this inequality in \eqref{4} we obtain $2 H(2\widetilde{T}) \ge H(t)$. Since $H''(t)=h'(t)\ge0$ we have also, integrating two times on $[2\widetilde{T},t]$, that $H(t) \ge H(2\widetilde{T}) + H'(2\widetilde{T}) \{t-2\widetilde{T}\}$. These two relations give us the estimate
	\begin{equation*}
	t \le 2\widetilde{T} + {H(2\widetilde{T})}/{H'(2\widetilde{T})} \quad\text{for $t>2\widetilde{T}$.}
	\end{equation*}
\par
	Finally, observe that
	\begin{equation*}
	\left({H(t)}/{H'(t)}\right)' = 1-{H(t)H''(t)}/{(H'(t))^2} \le 1
	\end{equation*}
	from which, integrating on $[\widetilde{T},2\widetilde{T}]$, we get
	\begin{equation*}
	{H(2\widetilde{T})}/{H'(2\widetilde{T})} \le  {H(\widetilde{T})}/{H'(\widetilde{T})} + 2\widetilde{T} - \widetilde{T}
	= \widetilde{T},
	\end{equation*}
	and so we have $t \le 3\widetilde{T}$.
	Therefore the proof of the lemma is completed.
\end{proof}

%%%%%%%%

%%%%%%%%%%%%%%%%%%%
%%%%%% section 4 %%%%%%%
%%%%%%%%%%%%%%%%%%%

%%%%%%%%%%%%%%%%%%%%%%%%%%%%
%%%%%%%%%%%%%%%%%%%%
%%%%%%%%%%%%%%%%%%%%%%%%%%%%

\section{Proof of Theorem \ref{thm:alpha<1}}
Following the idea in \cite{LST} and \cite{LT}, we introduce the multiplier
\begin{equation*}
\label{eq:m}
m(t):= \exp\left(\mu_1 \frac{(1+t)^{1-\beta}}{1-\beta}\right).
\end{equation*}
Clearly, $1 \geq m(t) \geq m(0) >0$ for $t \geq 0$.
Let us define the functional
\begin{equation*}
F_0(t) := \int_{\R^n} u(x,t)dx,
\end{equation*}
and then
$
F_0(0) = \e\int_{\R^n} f(x)dx,
F'_0(0) = \e\int_{\R^n}g(x)dx
$
are non-negative and do not both equal to zero, due to the hypothesis for the initial data.

Taking into account of (\ref{support}), we choose the test function $\phi=\phi(x,s)$ in the definition of energy solution \eqref{eq:energy_solution} such that it satisfies
$
\phi\equiv 1$ in $\{(x,s)\in \R^n\times[0,t]:|x|\le s+R\},
$
to get
\[
\begin{split}
&\int_{\R^n}u_t(x,t)dx-\int_{\R^n}u_t(x,0)dx+\int_0^tds\int_{\R^n}\frac{\mu_1}{(1+s)^\beta}u_t(x,s)dx \\
=& \int_0^t \int_{\R^n} \frac{\mu_2}{(1+s)^{\alpha+1}} u(x,s) dx
+\int_0^tds\int_{\R^n}|u(x,s)|^pdx,
\end{split}
\]
which yields, by taking derivative with respect to $t$,
\begin{equation}
\label{2}
F_0''(t)+\frac{\mu_1}{(1+t)^\beta}F_0'(t)
= \frac{\mu_2}{(1+t)^{\alpha+1}}F_0(t) + \int_{\R^n}|u(x,t)|^pdx.
\end{equation}
Multiplying both sides of \eqref{2} with $m(t)$ yields
\begin{equation}
\label{3}
\left\{m(t)F'_0(t)\right\}'
= m(t)\frac{\mu_2}{(1+t)^{\alpha+1}} \, F_0(t) + m(t)\int_{\R^n}|u(x,t)|^pdx.
\end{equation}
Integrating the previous equation twice on $[0,t]$, we obtain
\begin{comment}
	\begin{equation}
	\label{eq:estF'}
	\begin{split}
	F'_0(t) \ge &\ m(0)F'_0(0) + m(0)\mu_2 \int_{0}^{t} \frac{F_0(s)}{(1+s)^{\alpha+1}} ds \\
	&+ m(0) \int_{0}^{t} ds \int_{\R^n}|u(x,s)|^pdx
	\quad\text{for $t\ge0$},
	\end{split}
	\end{equation}
	and integrating again we reach
\end{comment}
\begin{equation}
\label{eq:estF}
\begin{split}
F_0(t)
\ge
&\ F_0(0) + F'_0(0) \int_{0}^{t} \frac{ds}{m(s)}
+ \mu_2 \int_{0}^{t} \frac{ds}{m(s)} \int_{0}^{s} \frac{m(\s)F_0(\s)}{(1+\s)^{\alpha+1}} d\s \\
&+ \int_{0}^{t} \frac{ds}{m(s)} \int_{0}^{s} m(\s)d\s \int_{\R^n}|u(x,\s)|^pdx
\quad\text{for $t\ge0$}.
\end{split}
\end{equation}

By a comparison argument, we observe $F_0(t)>0$ for $t>0$, and consequently also $F_0'(t) >0$ for $t>0$ by an integration of \eqref{3}. In fact, if $F_0(0)>0$, then $F_0$ is strictly positive for at least small times. Supposing that $t_0$ is the smallest zero time of $F_0$, calculating \eqref{eq:estF} in $t_0$ we get a contradiction. If $F_0(0)=0$, then $F'_0(0)>0$ and so $F_0'(t)>0$ for at least small time; due to the fact that $F_0$ is strictly increasing we then conclude that it is  positive for at least small time $t>0$. Supposing that $t_0>0$ is the smallest zero point of $F_0$, calculating \eqref{eq:estF} in $t_0$ we get again a contradiction.

Moreover observe that if
$F_0'(0)=0$, neglecting the last term on the right-hand side of \eqref{eq:estF}, and noting that $F$ is increasing and $m$ is bounded, we have
\begin{equation*}
F_0(t) \ge F_0(0) + 2^{-1}m(0)\mu_2 F_0(0) (1+t)^{-\max\{0,\alpha+1\}} t^2
\end{equation*}
and so $F_0(\widetilde{t}) \ge 2 F_0(0)$, if we choose $\widetilde{t}=\widetilde{t}(\mu_1,\mu_2,\alpha,\beta) >0$ such that
\begin{equation*}
2^{-1}m(0)\mu_2 (1+\widetilde{t})^{-\max\{0,\alpha+1\}} \widetilde{t}^2 = 1.
\end{equation*}

Now we need estimates for $\{mF_0'\}'$ and $F_0$. Neglecting the first term on the right-hand side of \eqref{3} and applying H\"older's inequality for the last term, there exists $C_1=C_1(n,p,R)>0$ such that, for $t\ge0$,
\begin{equation}
\label{F''}
	\left\{m(t)F'_0(t)\right\}'
	\ge m(0) C_1 (1+t)^{-n(p-1)} |F_0(t)|^p
	=: B(t)|F_0(t)|^p
	%\quad\text{for $t\ge0$}.
\end{equation}

Fix $t_0>0$ to be chosen later and consider the auxiliary function
\begin{equation*}
%\begin{split}
J(t) = \e J_0 + \e J_1 (t-t_0)
+ m(0)\mu_2 \int_{t_0}^{t} ds \int_{t_0}^{s}  \frac{J(\s)}{(1+\s)^{\alpha+1}} d\s
\quad\text{for $t \ge t_0$,}
%\end{split}
\end{equation*}
where $J_0 := \frac{1}{2}{\norm{f}_{L^1(\R^n)}}$, $J_1 := \frac{1}{2}{m(0)\norm{g}_{L^1(\R^n)}}$. By the similar way as above, we get by comparison argument that $F_0(t) \ge J(t)$ for $t\ge t_0$. Setting for the simplicity
$c:= m(0)\mu_2$, $q:=1-\alpha$,
the function $J$ satisfies
\begin{comment}
\begin{equation*}
\left\{
\begin{aligned}
& J''(t) = c (1+t)^{q-2} J(t)
\quad\text{for $t\ge t_0$}\\
& J(t_0) =\frac{\e}{2}\norm{f}_{L^1(\R^n)} =: \e J_0,
\quad
J'(t_0)= \frac{\e}{2} m(0)\norm{g}_{L^1(\R^n)} =: \e J_1.
\end{aligned}
\right.
\end{equation*}
\end{comment}
\begin{equation*}
	J''(t) = c (1+t)^{q-2} J(t)
	\quad\text{for $t\ge t_0$},
\end{equation*}
with $J(t_0) = \e J_0$, $J'(t_0)= \e J_1$.
One can check that the solution of this ordinary differential equation is
\begin{equation*}
J(t) = \e c_+ J_+(t) + c_- J_-(t)
\end{equation*}
where, setting $B^+_{1/q}:=I_{1/q}$ and $B^-_{1/q}:=K_{1/q}$ the modified Bessel functions respectively of the first and second kind with order $1/q$, we have
\begin{gather*}
J_{\pm}(t):= (1+t)^{1/2} B^{\pm}_{1/q} \left(2\frac{\sqrt{c}}{|q|} (1+t)^{q/2}\right),
\\
\begin{aligned}
c_\pm := &\ \pm \frac{2}{q} (1+t_0)^{-1/2}
\left[ (1+t_0) J_1 - J_0 \right]
B^{\mp}_{1/q}  \left(\frac{2\sqrt{c}}{|q|} (1+t_0)^{q/2}\right)
\\ &+ J_0 \frac{2\sqrt{c}}{|q|} (1+t_0)^{(q-1)/2}
B^{\mp}_{1+1/q} \left(\frac{2\sqrt{c}}{|q|}(1+t_0)^{q/2}\right).
\end{aligned}
\end{gather*}

Observe that $c_+>0$ at least for $t_0 > 0$ (independent of $\e$) large enough.
From the asymptotic expansions of the modified Bessel functions (see Section 9.7~in~\cite{AS}), when $t>0$ is large we have that
\begin{gather*}
J_\pm(t) = \frac{\sqrt{\pi^{\mp1}q}}{2c^{1/4}} (1+t)^{\frac{1}{2}-\frac{q}{4}} \exp\left(\pm2\frac{\sqrt{c}}{q} (1+t)^{\frac{q}{2}}\right)
%\times
\left[ 1 + O\left((1+t)^{-\frac{q}{2}}\right) \right]
.
\end{gather*}
Consequently, we can find constants $C_2, T_1 >0$ independent of $\e$, such that, for every $t \ge T_1$, the following estimate holds:
\begin{equation}
\label{eq:estG_alpha<1}
F_0(t) \ge \e C_2
(1+t)^{\frac{1+\alpha}{4}} \exp\left(2\frac{\sqrt{m(0)\mu_2}}{1-\alpha} (1+t)^{\frac{1-\alpha}{2}}\right) =: A(t) .
\end{equation}

%%%%%%%%%%%%%%%%%%%%%%%%%%%%%%%
%%%%%%%%%%%%%%%%%%%%%%%%%%%%%%%

Thanks to estimates \eqref{F''} and \eqref{eq:estG_alpha<1}, we can apply Lemma~\ref{lem:kato}.
\begin{comment}
whit
\begin{gather*}
A(t) = \e C_3
(1+t)^{(1+\alpha)/4} \exp\left(2\frac{\sqrt{m(0)\mu_2}}{1-\alpha} (1+t)^{(1-\alpha)/2}\right), \\
%
B(t) = m(0)C_1 (1+t)^{-n(p-1)}.
\end{gather*}
\end{comment}
Fix $\delta :=(p-1)/4$
\begin{comment}
and define
\begin{equation*}
h(t) := B(t)^{1/2} A(t)^{(p-1)/2-\delta},
\quad
h_0(t) := B(t)^{1/2} A(t)^{(p-1)/2}.
\end{equation*}
Observe that
\end{comment}
and, using the definition \eqref{def:h} of $h$, observe that
\begin{equation*}
\begin{split}
h'(t)
%&= B^{1/2}(t) A^{(p-1)/2-\delta}(t) \left[\frac{B'(t)}{2B(t)} + \left(\frac{p-1}{2}-\delta \right) \frac{A'(t)}{A(t)}\right] \\
&= B^{1/2}(t) A^{(p-1)/2-\delta}(t) (1+t)^{-1} \\
&\times \left\{ -\frac{n(p-1)}{2}
+ \left( \frac{p-1}{2}-\delta \right)
\left[\frac{1+\alpha}{4} + \sqrt{m(0)\mu_2} (1+t)^{\frac{1-\alpha}{2}}
\right]
\right\},
\end{split}
\end{equation*}
so we can find a time $T_2=T_2(n,p,\mu_1,\mu_2,\alpha,\beta) \ge 0$ such that $h'(t) > 0$ for $t \ge T_2$. Then we can choose $\widetilde{T}_0 = \max\{ T_1, T_2\}$.

Let us set $\widetilde{T} \equiv \zeta-1$ and
$C = C_3 [{\delta m(0) \sqrt{C_1}}/{(2\sqrt{p+1})}]^{2/(p-1)},$
\begin{comment}
\begin{equation*}
C = C_3 \left(\frac{\delta m(0) \sqrt{C_1}}{2\sqrt{p+1}}\right)^{2/(p-1)},
\end{equation*}
\end{comment}
where $\zeta\equiv \zeta(\overline{\e})$ is the larger solution to \eqref{eq:defb} with $\overline{\e}=C\e$. There exists $\e_0 >0$ such that, for $0< \e \le \e_0$, we have $\widetilde{T} \equiv \zeta-1 \ge \max\{\widetilde{T}_0,\widetilde{T}_1\}$, where $\widetilde{T}_1$, independent of $\e$, is defined as in the statement of the Lemma. We can also suppose $\widetilde{T}\ge1$, so that $\zeta-1 \ge \zeta/2$. Therefore, one can check that \eqref{eq:T1estimate} holds
\begin{comment}
\begin{equation*}
\begin{split}
A(\widetilde{T})^\delta \widetilde{T} \, h(\widetilde{T})
&\ge
\frac{1}{2} C_3^{\frac{p-1}{2}} \sqrt{m(0)C_1} \e^{\frac{p-1}{2}}
\zeta^{1-\frac{n(p-1)}{2}+\frac{1+\alpha}{4}\cdot\frac{p-1}{2}} \exp\left(\frac{2\sqrt{m(0)\mu_2}}{1-\alpha} \cdot\frac{p-1}{2} \zeta^{\frac{1-\alpha}{2}}\right) \\
%&=
\frac{1}{\delta}\sqrt{\frac{p+1}{m(0)}} \left[ C\e \zeta^{\frac{2}{p-1}-n+\frac{1+\alpha}{4}} \exp\left(\frac{2\sqrt{m(0)\mu_2}}{1-\alpha} \zeta^{(1-\alpha)/2}\right) \right]^{\frac{p-1}{2}} \\
&= \frac{1}{\delta}\sqrt{\frac{p+1}{m(0)}},
\end{split}
\end{equation*}
\end{comment}
and so the maximal existence time $T$ of $F_0$ satisfies $T \le 3\widetilde{T} \le 3\zeta$.
The proof of the Theorem~\ref{thm:alpha<1} is completed.
%%%%%%%%%%%%%%%%%%%%%%%%%%%%%%%%%
% Appendix
%%%%%%%%%%%%%%%%%%%%%%%%%%%%%%%%%

\section{Appendix}

In this section we are going to show the local existence and finite speed of propagation property for energy solution to our problem, as stated in the second point of Remark \ref{rem1}. In the following, the positive constant $C$ may vary from line to line.
We assume that $p\le n/(n-2)$ when $n\ge3$.

Let us denote the function space
\begin{equation*}\label{functionspace}
\begin{aligned}
B_{T, K} :=
\big\{ & \phi \in C\left([0, T), H^1(\R^n)\right)\cap C^1\left([0, T), L^2(\R^n)\right) \colon \\
& \supp \phi \subset \{(x,t)\in\R^n\times[0,T)\colon |x|\le t+R \} ,
\,\,
 \|\phi\|_{B_{T, K}}\le K \big\},
\end{aligned}
\end{equation*}
where $T, R, K$ are fixed positive constants and
\begin{equation*}\label{energynorm}
\|\phi\|_{B_{T, K}} := \sup_{t\in[0, T)}E_{\phi}^{1/2}(t),
\quad
E_{\phi}(t):=\frac{1}{2}\int_{\R^n}(\phi_t^2+|\nabla \phi|^2)dx.
\end{equation*}
It can be proved that $B_{T, K}$ is a Banach space.

Consider the following Cauchy problem for $v\in B_{T,K}$
\begin{equation}
\label{Cau_map}
\left\{
\begin{aligned}
& u_{tt} - \Delta u=|v|^p+m^2(t)v-b(t)v_t =:F_v(x, t),~ \text{in $\R^n\times(0,T)$}, \\
& u(x,0)=\e f(x), \quad u_t(x,0)=\e g(x), \quad x\in\R^n,
\end{aligned}
\right.
\end{equation}
where we set for the simplicity
\[
m^2(t)=\frac{\mu_2}{(1+t)^{\alpha+1}},
\quad
b(t)=\frac{\mu_1}{(1+t)^\beta}.
\]
However, all the calculations below are trivially generalized for any $m^2,b \in C([0,T))$.

We want to show that the map
\[
M \colon v \mapsto u=Mv,~~v\in B_{T,K},
\]
is a contraction. Note that, for $v\in B_{T,K}$, by Gagliardo-Nirenberg inequality and Poincar\'e inequality, we have
\[
\|v\|_{L^{2p}(\R^n)}\le C\|v\|_{L^2(\R^n)}^{1-\theta(2p)}\|\nabla v\|_{L^2(\R^n)}^{\theta(2p)},
\quad\theta(2p):=n\left(\frac{1}{2}-\frac{1}{2p}\right),
\]
for $p\le n/(n-2)$ when $n\ge3$, and
\[
\|v\|_{L^2(\R^n)}\le C(t+R)\|\nabla v\|_{L^2(\R^n)},
\]
which imply
\begin{equation}\label{GN+P}
\|{v}\|_{L^{2p}(\R^n)}
\le C (t+R)^{1-\theta(2p)} \|\nabla {v}\|_{L^2(\R^n)}
\le C (t+R)^{1-\theta(2p)} E_{{v}}^{1/2}.
\end{equation}
In particular, for fixed $T>0$, we can check that
\[
F_v(x, t)\in L^2\left(\R^n\times [0, T)\right).
\]

Let us start proving that the map $M$ is onto. Firstly, we show the finite speed propagation of the energy solution, i.e.
\[
\supp u \subset \{(x,t)\in \R^n\times[0,T) \colon |x|\le t+R\},
\]
by using the density argument similarly to \cite{Struwe}. By density of $C_0^{\infty}(\R^n)$ in $L^2(\R^n)$, we can approximate the energy data $f, g$ by sequences of smooth and compactly supported functions
$\{f_m\}_{m\in\N}, \{g_m\}_{m\in\N}$ in the energy norm $H^1(\R^n)$ and $L^2(\R^n)$ respectively. Noting that $F_v(x, t)$ has compact support, we can find also a sequence of smooth and compactly supported functions $\{F_{v,m}\}_{m\in\N}$ converging to $F_v$ in the norm $L^2\left(\R^n\times [0, T)\right)$. Let $u_m$ be the smooth solution of the problem
\begin{equation}\label{Cau_approx}
\left\{
\begin{array}{l}
(u_m)_{tt}-\Delta u_m=F_{v,m}(x,t) \quad \text{in $\R^n\times(0,T)$}\\
u(x,0)=\e f_m(x), \quad u_t(x,0)= \e g_m(x)
\quad \text{in $\R^n$.}
\end{array}
\right.
\end{equation}
Fix $(x_0, t_0)\in \R^n\times (0, T)$ with $|x_0|\ge t_0+R$ and set
\begin{equation*}
	C_{(x_0,t_0)} := \{ (x,t) \in \R^n\times[0,T) \colon |x-x_0|\le t_0-t \},
\end{equation*}
the backward cone with vertex in $(x_0,t_0)$.
Then, denote the energy on a time-section of the cone as
\begin{equation}\label{local_ennorm}
E_{t_0-t}(t,u(t)):=\frac{1}{2}\int_{B_{t_0-t}(x_0)}(u_t^2+|\nabla u|^2)dx,
\end{equation}
where $B_r(x_0):=\{x\in\R^n:|x-x_0|\le r\}$.
The standard space-time divergence form yields
a local energy inequality
\begin{equation*}
E^{1/2}_{t_0-t}(t,u_m(t)) \le E^{1/2}_{t_0}(0,u_m(0)) + \int_{0}^{t} \|F_{v,m}(\cdot,s)\|_{L^2(B_{t_0-s}(x_0))} \,ds.
\end{equation*}
Applying the previous inequality to the difference $u_m-u_n$ of two solutions of \eqref{Cau_approx}, we have that $\{u_m(\cdot,t)\}_{m\in\N}$ is a Cauchy sequence in the norm \eqref{local_ennorm} uniformly in $t\in[0,t_0]$. Hence the limit $u$ is an energy solution satisfying
\[
E^{1/2}_{t_0-t}(t,u(t))
\le E^{1/2}_{t_0}(0,u(0))
+\int_0^t \|F_v(\cdot,s)\|_{L^2(B_{t_0-t}(x_0))} \,ds,
\]
which gives us the fact that
\[
f(x)\equiv g(x)\equiv0\quad\mbox{in}\ C_{x_0,t_0}\cap\{t=0\}
\quad\mbox{and}\quad
F_v\equiv 0\ \mbox{in}\ C_{(x_0,t_0)}
\]
and  Poincar\'e inequality imply
\[
u\equiv0\ \mbox{in}\ C_{(x_0,t_0)}.
\]

Next, we show that $\|Mv\|_{B_{T,K}} \le K$. It is easy to obtain
\begin{equation}
\label{multiply}
\frac{\p}{\p t}\left(\frac{1}{2}(u_t^2+|\nabla u|^2)\right)
=\mbox{div}(u_t\nabla u)+|v|^pu_t+m^2(t)vu_t-b(t)v_tu_t.
\end{equation}
Exploiting \eqref{GN+P} we get the estimate
\begin{align*}
\int_{\R^n}|v|^p|u_t|dx
&\le
\left(\int_{\R^n}|v|^{2p}dx\right)^{1/2}\sqrt{2}E^{1/2}_{u}(t)\\
%&\le \sqrt{2} \|v\|_{L^{2p}(\R^n)}^p E^{1/2}_{u}(t) \\
&\le
C(t+R)^{{p}\{1-\theta(2p)\}}E^{p/2}_{v}(t)E^{1/2}_{u}(t).
\end{align*}
Moreover, it is trivial that
\[
\int_{\R^n}|v||u_t|dx\le C(t+R)E^{1/2}_{v}(t)E^{1/2}_{u}(t)
\]
and
\[
\int_{\R^n}|v_t||u_t|dx\le2E^{1/2}_{v}(t)E^{1/2}_{u}(t).
\]

Integrating \eqref{multiply} over $\R^n \times [0,t]$ and using the divergence theorem and the estimates above, we obtain
\begin{equation*}
E_{u}(t)\le
E_{u}(0)
+C\int_0^t a_K(s)
E_{u}(s)^{1/2}ds,
\end{equation*}
where
\begin{equation*}
	a_K(t):=  K^{p} (t+R)^{p\{1-\theta(2p)\}}+ K (t+R) m^2(t)+ K b(t) ,
\end{equation*}
which yields, by Bihari's inequality, that for some positive constant $\gamma$
\begin{equation*}
\begin{aligned}
E_{u}(t)^{1/2} &\le E^{1/2}_{u}(0)
+C\int_0^t a_K(s) ds\\
&\le E_u^{1/2}(0)+C \max\{K, K^{p}\} T(1+T)^{\gamma}.
\end{aligned}
\end{equation*}
Hence we can choose $K$ large enough and $T$ small enough such that
\[
E_u^{1/2}(0)\le \frac{K}{2}
\quad
\text{and}
\quad
C \max\{K, K^{p}\} T(1+T)^{\gamma}\le \frac{K}{2},
\]
and then $E_u^{1/2}(t)\le K$.

\par

%%%%%

Finally, we can prove the contraction of the map $M$ in a similar way. Fixed $v_1, v_2\in B_{T, K}$, let
\begin{equation*}
u_1=Mv_1, \quad u_2=Mv_2
\quad
\text{and}
\quad
\overline{u}=u_1-u_2,
\quad
\overline{v}=v_1-v_2.
\end{equation*}
We have that $\overline{u}$ satisfies the problem
\begin{equation*}
\left\{
\begin{aligned}
& \overline{u}_{tt}-\Delta\overline{u}=|v_1|^p-|v_2|^p+m^2(t)\overline{v}-b(t)\overline{v}_t,
 \quad\text{in $\R^n\times(0,T)$}, \\
& u(x,0)= u_t(x,0) \equiv 0, \quad x\in\R^n,
\end{aligned}
\right.
\end{equation*}
and the equation
\begin{multline}\label{multiply1}
\frac{\p}{\p t}\left(\frac{1}{2}(\overline{u}_t^2+|\nabla \overline{u}|^2)\right) \\
=
\text{div}(\overline{u}_t\nabla \overline{u})
+\left(|v_1|^p-|v_2|^p\right)\overline{u}_t
+m^2(t)\overline{v}\overline{u}_t-b(t)\overline{v}_t\overline{u}_t.
\end{multline}
Observe that, by \eqref{GN+P}, it holds
\begin{equation*}\label{nonlin}
\begin{aligned}
\int_{\R^n} ||v_1|^p-|v_2|^p|\,|\overline{u}_t| \,dx
&\le
C\int_{\R^n}|v_1-v_2|(|v_1|+|v_2|)^{p-1}|\overline{u}_t| dx\\
&\le C \| |\overline{v}|(|v_1|+|v_2|)^{p-1}\|_{L^2}\|\overline{u}_t\|_{L^2}\\
&\le C
\|\overline{v}\|_{L^{2p}}
\left(\| v_1\|_{L^{2p}} + \|v_2\|_{L^{2p}}\right)^{p-1}
\|\overline{u}_t\|_{L^2},\\
&
\le C K^{p-1} (t+R)^{p\{1-\theta(2p)\}}
E^{1/2}_{\overline{v}} E^{1/2}_{\overline{u}}.
\end{aligned}
\end{equation*}
Then, integrating \eqref{multiply1} on $\R^n\times[0,t]$, exploiting again the Bihari's inequality and proceeding similarly as above, we reach the estimate
\begin{equation*}
	\|\overline{u}\|_{B_{T,K}} \le
	C
	\max\{1,K^{p-1}\} T(1+T)^\gamma \|\overline{v}\|_{B_{T,K}},
\end{equation*}
from which, choosing $T$ small enough, we infer that $M$ is a contraction.
The proof of the desired local existence is now completed.

\section*{Acknowledgement}
The first author is partially supported by Natural Science Foundation of Zhejiang Province(LY18A010008), NSFC(11771194), high level talent project of Lishui City (2016RC25), the Scientific Research Foundation of the First-Class Discipline of Zhejiang Province(B)(201601).
The third author is partially supported by the Grant-in-Aid for Scientific Research (B)(No.18H01132), Japan Society for the Promotion of Science.

%%%%%%%%%%%%%%%%%%%%%%%%%%%%%%%%%
% References
%%%%%%%%%%%%%%%%%%%%%%%%%%%%%%%%%


\begin{thebibliography}{99}

%%%%%%%%%%%%%%%%%%%%%%%%%%%%%%%%%%%%%%%
%\bibitem{B1}
%T. Bridgeland,
%Equivalences of triangulated categories and Fourier-Mukai transforms,
%Bull. London Math. Soc., {\bf 31} (1999), 25--34.
%%%%%%%%%%%%%%%%%%%%%%%%%%%%%%%%%%%%%%%%
%\bibitem{Sa}
%I. Satake,
%Algebraic Structures of Symmetric Domains,
%Publ. Math. Soc. Japan, {\bf 14},
%Iwanami Shoten, Tokyo;
%Princeton University Press, Princeton, N.J., 1980.
%%%%%%%%%%%%%%%%%%%%%%%%%%%%%%%%%%%%%%%%%
%\bibitem{Vo} D. A. Vogan,  Jr.,
%A Langlands classification for unitary representations,
%In: Analysis on Homogeneous Spaces and Representation Theory of Lie Groups,
%Okayama-Kyoto,  1997,
%(eds. T. Kobayashi, M. Kashiwara, T. Matsuki, K. Nishiyama and T. Oshima),
%Adv. Stud. Pure Math., {\bf 26}, Math. Soc. Japan, 2000, pp. 299--324.
%%%%%%%%%%%%%%%%%%%%%%%%%%%%%%%%%%%%%%%%%%

\bibitem{AS}
Abramowitz, M., Stegun, I.A.:
Handbook of mathematical functions with formulas, graphs, and mathematical table.
Vol. 2172. Courier Corporation, 1965

%\bibitem{DABI} D'Abbicco, M.: The threshold of effective damping for semilinear wave equations. Math. Methods Appl. Sci. {\bf 38}, 1032-1045 (2015)

%\bibitem{DL1} D'Abbicco, M., Lucente, S.: A modified test function method for damped wave equations. Adv. Nonlinear Stud. {\bf 13}, 867-892 (2013)

%\bibitem{DLR13} D'Abbicco, M., Lucente, S., Reissig, M.: Semi-linear wave equations with effective damping. Chin. Ann. Math. Ser. B {\bf 34}, 345-380 (2013)

\bibitem{DLR15}
D'Abbicco, M., Lucente, S., Reissig, M.:
A shift in the Strauss exponent for semilinear wave equations with a not effective damping.
J. Differ. Equ. {\bf 259}, 5040-5073 (2015)

\bibitem{Fujita}
H. Fujita:
On the blowing up of solutions of the Cauchy problem for $u_t=\Delta u+u^{1+\alpha}$, J. Fac.
Sci. Univ. Tokyo Sect. I, 13, 109-124 (1966)

\bibitem{FIW}
Fujiwara, K., Ikeda, M., Wakasugi, Y.:
Estimates of lifespan and blow-up rate for the wave equation with a time-dependent damping and a power-type nonlinearity. arXiv:1609.01035 (2016)

\bibitem{II} Ikeda, M., Inui, T.:
The sharp estimate of the lifespan for the semilinear wave equation with time-dependent damping.
Differ. Integral Equ. 32.1/2, 1-36 (2019)

%\bibitem{IO} Ikeda, M., Ogawa, T.: Lifespan of solutions to the damped wave equation with a critical nonlinearity. J. Differ. Equ. {\bf 261}, 1880-1903 (2016)

%\bibitem{IS} Ikeda, M., Sobajima, M.: Life-span of solutions to semilinear wave equation with time-dependent critical damping for specially localized initial data. Math. Ann., (2018), https://doi.org/10.1007/s00208-018-1664-1.

\bibitem{IWnew}
Ikeda, M., Wakasugi, Y.:
Global well-posedness for the semilinear damped wave equation with time dependent damping in the overdamping case.
arXiv:1708.08044 (2017)

\bibitem{LST}
Lai, N.-A., Schiavone, N.M., Takamura, H.:
Wave-like blow-up for semilinear wave equations with scattering damping and negative mass term,
D'Abbicco, M., Ebert, M., Georgiev, V., and Ozawa, T. ed.,
\lq\lq New Tools for Nonlinear PDEs and Application",
Trends in Mathematics, 217-240, Birkh\"auser (2019)

\bibitem{LTW}
Lai, N.-A., Takamura, H., Wakasa, K.:
Blow-up for semilinear wave equations with the scale invariant damping and super-Fujita exponent.
J. Differ. Equ. {\bf 263}(9), 5377-5394 (2017)
%\url{http://dx.doi.org/10.1016/j.jde.2017.06.017}.

\bibitem{LT}
Lai, N.-A., Takamura, H.:
Blow-up for semilinear damped wave equations with subcritical exponent in the scattering case.
Nonlinear Anal. {\bf 168}, 222-237 (2018)

%\bibitem{LNZ12} Lin, J., Nishihara, K., Zhai, J.: Critical exponent for the semilinear wave equation with time-dependent damping. Discrete Contin. Dyn. Syst. Ser. A {\bf 32}, 4307-4320 (2012)

%\bibitem{NPR} Nunes do Nascimento, W., Palmieri, A., Reissig, M.: Semi-linear wave models with power non-linearity and scale-invariant time-dependent mass and dissipation. Math. Nachr. (2016) doi: 10.1002/mana.201600069

%\bibitem{N1} Nishihara, K.: Asymptotic behaviour of solutions to the semilinear wave equation with time-dependent damping. Tokyo J. Math. {\bf 34}, 327-343 (2011)

\bibitem{NPR}
Nascimento, N., Palmieri, A., Reissig, M.:
Semi-linear wave models with power non-linearity and scale-invariant time-dependent mass and dissipation.
Math. Nachr. {\bf 290}(11-12), 1779-1805 (2017)
%doi: 10.1002/mana.201600069	

\bibitem{Pal} Palmieri, A.:
Global existence of solutions for semi-linear wave equation with scale-invariant damping and mass in exponentially weighted spaces. J. Math. Anal. Appl. {\bf 461}(2), 1215-1240 (2018)

\bibitem{P-R}
Palmieri, A., Reissig, M.:
A competition between Fujita and Strauss type exponents for blow-up of semi-linear wave equations with scale-invariant damping and mass.
J. Differ. Equ., {\bf 266}(2-3), 1176-1220 (2019)

\bibitem{Strauss}
Strauss, W. A.:
Nonlinear scattering theory at low energy.
J. Funct. Anal. {\bf 41}(1), 110-133 (1981)

\bibitem{Struwe}
Struwe, M.:
Semilinear wave equations.
Bulletin of the American Mathematical Society {\bf 26}(1), 53-85 (1992)

\bibitem{T}
Takamura, H.:
Improved Kato's lemma on ordinary differential inequality and its application to semilinear wave equations.
Nonlinear Anal. Theory Methods Appl. {\bf 125}, 227-240 (2015)

\bibitem{TL1} Tu, Z., Lin, J.:
A note on the blowup of scale invariant damping wave equation with sub-Strauss exponent.
arXiv:1709.00866 (2017)

%\bibitem{TL2} Tu, Z., Lin, J.: Life-span of semilinear wave equations with scale-invariant damping: critical Strauss exponent case. arXiv:1711.00223.

\bibitem{WY} Wakasa, K., Yordanov, B.:
On the nonexistence of global solutions for critical semilinear wave equations
with damping in the scattering case.
Nonlinear Anal. TMA {\bf 180}, 67-74 (2019)

\bibitem{WY14_scale} Wakasugi, Y.:
Critical exponent for the semilinear wave equation with scale invariant damping.
Fourier analysis, Birkh\"{a}user, Cham, 375-390 (2014)

%\bibitem{WY17} Wakasugi, Y.: Scaling variables and asymptotic profiles for the semilinear damped wave equation with variable coefficients. J. Math. Anal. Appl. {\bf 447}, 452-487 (2017)

\bibitem{Wang}
Wang, C., Liu, M.:
Global existence for semilinear damped wave equations in relation with the Strauss conjecture.
arXiv:1807.05908 (2018)

\bibitem{Wir1}
Wirth, J.:
Solution representations for a wave equation with weak dissipation.
Math. Methods Appl. Sci. {\bf 27}, 101-124 (2004)

\bibitem{Wir2}
Wirth, J.:
Wave equations with time-dependent dissipation. I. Non-effective dissipation.
J. Differ. Equ. {\bf 222}, 487-514 (2006)

\bibitem{Wir3} Wirth, J.:
Wave equations with time-dependent dissipation. II. Effective dissipation.
J. Differ. Equ. {\bf 232}, 74-103 (2007)

\bibitem{YG}
Yagdjian, K., Galstian, A.:
The Klein-Gordon equation in anti-de Sitter spacetime.
Rend. Semin. Mat. Univ. Politec. Torino {\bf 67}(2), 271-292 (2009)

%\bibitem{YZ06} Yordanov, B., Zhang, Q.S.: Finite time blow up for critical wave equations in high dimensions. J. Funct. Anal. {\bf 231}, 361-374 (2006)

%\bibitem{Z} Zhang, Q.S.: A blow-up result for a nonlinear wave equation with damping: the critical case. C. R. Math. Acad. Sci. Paris, S\'er. I {\bf 333}, 109-114 (2001)

\end{thebibliography}
\end{document}